\RequirePackage{fix-cm}
\documentclass[smallextended]{svjour3}       
\smartqed  
\usepackage{graphicx}
\usepackage{amsmath}
\usepackage{amssymb}
\begin{document}

\title{Optimal Control Problems for Evolutionary Variational
Inequalities with Volterra type Operators
}
\subtitle{Optimal Control Problems for Variational Inequalities}


\author{Mykola Bokalo         \and
        Olha Sus 
}


\institute{Mykola Bokalo \at
              Department of Mechanics and Mathematics \\
              Ivan Franko National University of Lviv \\
              Lviv, Ukraine \\
              \email{mm.bokalo@gmail.com}           
           \and
           Olha Sus \at
           Department of Mathematics and Statistics \\
           University of Alaska Fairbanks \\
           Fairbanks, AK, USA \\
           \email{osus@alaska.edu}   
}

\date{Received: date / Accepted: date}

\maketitle

\begin{abstract}
In this paper, we consider the optimal control problem for a class of evolution inclusions with Volterra type operators, which can be history-dependent. We establish the existence of a solution to the stated optimal control problem under some hypothesis on data-in. The motivation for this work comes from the optimal control problems for variational inequalities arising in the study of injection molding processes, contact mechanics, principles of electro-wetting on dielectric, and others.
\end{abstract}

\section{Introduction}
\label{intro}
In this paper, we consider optimal control problems for parabolic variational inequalities (subdifferential inclusions) with Volterra type operators. Let us first introduce an example of the problem which is studied here.

Let $n\in \mathbb{N}$, $p\geq2$, $T>0$ be given numbers, $\Omega$ be a bounded domain in $\mathbb{R}^n$ with a  boundary $\partial \Omega$. We put  $Q:=\Omega\times (0,T)$, $\Sigma:=\partial\Omega\times (0,T)$,  $\Pi:=\{(t,s)\,|\, t\in (0,T), s\in (0,t)\}$.

Let $L^p(\Omega)$, $L^p(Q)$ be standard Lebesgue spaces. We denote by $W^{1,p}(\Omega)=\{v\in L^p(\Omega) \ | \ v_{x_i}\in L^p(\Omega),\,\,i=\overline{1,n}\}$ a standard Sobolev space with the norm
$\|v\|_{W^{1,p}(\Omega)}:=\Big(\sum_{i=1}^{n}\|v_{x_i}\|_{L^p(\Omega)}^p
+\,\|v\|_{L^p(\Omega)}^p\Big)^{1/p}.$

We consider operator $\widehat{\mathcal{B}}: C([0,T];L^2(\Omega))\to L^\infty(0,T;L^2(\Omega))$, defined by the following rule: for each function $z(x,t),\, (x,t)\in Q,$ from $C([0,T];L^2(\Omega))$ we have
\begin{equation*}
\widehat{\mathcal{B}}(z)(x,t):=
\int\limits_0^t \widehat{b}(x,t,s)\,z(x,s)\,ds,\, (x,t)\in Q,
\end{equation*}
where $\widehat{b}\in L^\infty(\Omega\times\Pi)$ is given.
Note, that for almost every (a.e.) $t\in (0,T)$ and every $z\in C([0,T];L^2(\Omega))$ we obtain
\begin{equation}\label{opVolterraexample1}
\|\widehat{\mathcal{B}}(z)(\cdot,t)\|_{L^2(\Omega)}=
\Big\|\int\limits_0^t \widehat{b}(\cdot,t,s)\,z(\cdot,s)\,ds\Big\|_{L^2(\Omega)}\leq
$$ $$\leq \int\limits_0^t \|\widehat{b}(\cdot,t,s)\,z(\cdot,s)\|_{L^2(\Omega)}\,ds \leq
\widehat{L} \int\limits_0^t \|z(\cdot,s)\|_{L^2(\Omega)}\,ds,
\end{equation}
where $\widehat{L}:=\mathop{\textrm{ess}\sup}\limits_{(x,t,s)\in \Omega\times\Pi}|\widehat{b}(x,t,s)|$.
Operator $\widehat{\mathcal{B}}$ is called a Volterra type operator.

Let $U:=L^2(Q)$ be a space of controls, and $U_{\partial}$ be a closed subset of $U$. For example, $U_{\partial}:=\{u\in U \ | \ m\leq u\leq M\ \textrm{a.e. on} \ Q\}$ for given  $m,M\in\mathbb{R}$. The set $U_{\partial}$ is called  the set of admissible controls. Also, let $K$ be a convex closed set in $W^{1,p}(\Omega)$ which contains $0$. For example, we set $K=\{v\in W^{1,p}(\Omega)\,|\, v\geq 0\ \textrm{almost every on} \ \Omega\}$.

For a given control $u\in U_{\partial}$ the state of evolutionary system  is described by a function
$y(x,t),\, (x,t)\in\overline{Q},$ ( also it can be denoted by
$y$, or $y(u)$, or $y(x,t;u),\, (x,t)\in \overline{Q}$) such that $y\in L^p(Q)\cap C([0,T];L^2(\Omega)),\, y_{x_i}\in L^p(Q)$, $i=\overline{1,n}$, $y_t\in L^2(Q)$, and
\begin{equation}\label{8I_pochumova}
y|_{t=0}=y_0(x),\,\quad x\in \Omega,
\end{equation}
and, for a.e. $t\in (0,T]$, $u(\cdot,t)\in K$ and
\begin{equation}\label{8I_varinequality1}
\int\limits_{\Omega_t}\big\{ y_t(v-y) + |\nabla y|^{p-2}\nabla y\nabla (v-y) +
 |y|^{p-2}y(v-y)+
 $$
 $$
 +(v-y)\int_0^t \widehat{b}(x,t,s)y(x,s)\,ds\big\}\,dx\geq\int\limits_{\Omega_t}(f+u)(v-y)\,dx\quad \forall\,v\in K,
\end{equation}
where  $f\in L^2(Q),\,y_0\in L^2(\Omega)$ are given, $\Omega_t:=\Omega\times \{t\}$,
$\nabla y:=(y_{x_1},\ldots,y_{x_n})$.

The cost functional $J:U\rightarrow\mathbb{R}$ is defined by the following rule
\begin{equation}\label{ast}
J(u)=\|y(\cdot,T;u)-y_T(\cdot)\|^2_{L^2(\Omega)}+\mu\|u\|^2_{L^2(Q)},
\end{equation}
where $\mu>0$, $y_T\in L^2(\Omega)$ are given.

The \textbf{optimal control problem} is to find a control $u_{\ast}\in U_{\partial}$ such that
\begin{equation}\label{Juast}
J(u_{\ast})=\inf_{u\in U_{\partial}}J(u).
\end{equation}
\quad
As it will be shown later, there exists a solution to the stated optimal control problem (\ref{Juast}).

We remark that formulated above problem can be written in more abstract way.
Indeed, after an appropriate identification of functions and functionals, one can write
\begin{equation*}
W^{1,p}(\Omega)\subset L^2(\Omega) \subset (W^{1,p}(\Omega))',
\end{equation*}
where all embeddings are continuous and dense, provided $W^{1,p}(\Omega)$ is reflexive, and the first embedding is compact. Clearly, for any
$h\in L^2(\Omega)$ and $v\in W^{1,p}(\Omega)$ we have
$\langle h,v\rangle=(h,v)$, where by $\langle \cdot,\cdot\rangle$ we denote the scalar product on the dual pair $\bigl[(W^{1,p}(\Omega))'$, $W^{1,p}(\Omega)\bigr]$, and  by $(\cdot,\cdot)$ we denote the scalar product in $L^2(\Omega)$.
Thus, we may use the notation $(\cdot,\cdot)$ instead of $\langle \cdot,\cdot\rangle$.

Now let us denote by $V:=W^{1,p}(\Omega)$, $H:=L^2(\Omega)$
and define  operator $A:V\to V'$  as follows
\begin{equation*}
(A(v),w):=\int\limits_\Omega\big[ |\nabla v|^{p-2}\nabla v\nabla w + |v|^{p-2}vw \big]\,dx, \quad v,w \in V.
\end{equation*}
Also, we use the following notations
\begin{equation*}
\widehat{b}(t,s):=\widehat{b}(\cdot,t,s), \ \ (t,s)\in \Pi,\quad f(t):=f(\cdot,t),\ \ t\in (0,T),
\end{equation*}
\begin{equation*}
u(t):=u(\cdot,t),\ \ y(t):=y(\cdot,t),\ \ y(t;u):=y(\cdot,t;u), \ \ t\in [0,T].
\end{equation*}

Then for a given control $u\in U_{\partial}$ the state of evolutionary system  is described by a function
$y \in L^p(0,T;V)$  such that $y' \in L^2(0,T;H)$, $y(0)=y_0$
and, for a.e.  $t\in (0,T)$,  $y(t) \in K$ and
\begin{equation}\label{8I_varinequality2}
  \bigl(y'(t)+A(y(t))+\int_0^t \widehat{b}(t,s)y(s)\,ds\,, v-y(t)\bigr)\geq\bigl(f(t)+u(t),v-y(t)\bigr) \quad \forall\, v\in K,
\end{equation}
where $f\in L^2(0,T;H),\ y_0\in H$ are given.

In this case (see (\ref{ast})), the cost functional $J:U\rightarrow\mathbb{R}$ is defined by
\begin{equation*}
J(u):=\|y(T;u)-y_T\|^2_H+\mu \|u\|^2_U,
\end{equation*}
where $y_T\in H$ is given.

We remark that variational inequality (\ref{8I_varinequality2}) can be written as a subdifferential  inclusion.

For this purpose we put $I_K(v):=0$ if $v\in K$, and $I_K(v):=+\infty$ if $v\in V\setminus K$, and also
\begin{equation*}
\Phi(v):=\frac{1}{p}\int\limits_\Omega \big[|\nabla v|^p+|v|^p\big]\,dx + I_K(v), \quad v\in V.
\end{equation*}
\qquad
It is easy to verify that the functional $\Phi:V\to \mathbb{R} \cup \{+\infty\}$ is convex and semi-lower-continuous.
By the known results (see., e.g., \cite[p. 83]{1.}) it follows that
the problem of finding a solution to the variational inequality (\ref{8I_varinequality2}) with initial condition $y(0)=y_0$
is equivalent to finding
{\it a function $y\in L^p(0,T;V)\cap C([0,T];H)$ such that $y'\in L^2(0,T;H)$, $y(0)=y_0$ and, for a.e. $t\in (0,T)$, $y(t) \in D(\partial\Phi)$ and
\begin{equation}\label{8I_varinequality3}
y'(t)+\partial \Phi(y(t))+\int_0^t b(t,s)y(s)\,ds \ni f(t)+u(t)\quad \text{in} \quad H,
\end{equation}
where $\partial\Phi:V\rightarrow2^{V'}$ is a subdifferential of the functional $\Phi$.}

The aim of the present paper is to solve the optimal control problem for inclusions type (\ref{8I_varinequality3}).

Optimal control problems for variational inequalities have been a subject of interest in the optimal control community starting from the 1980s. The motivation for this study comes from broad interesting applications. For instance, in \cite{2.}, one can find its application in the problem related to principles of electro-wetting on dielectric. Also, such problems appear at the simulation of various problems related to injection molding processes, contact mechanics, etc. (see, e.g.,\cite{3.}). Furthermore, it has its application in problems from economics, finance, optimization theory, and many others (see, e.g. \cite{4.}, \cite{5.}, \cite{6.}, and reference therein).

The optimal control of evolution problems has been extensively studied in the literature. In the book \cite{7.}, it was considered by Lions the optimal control of systems governed by partial differential equations. The existence and approximation of optimal solutions as well as the necessary optimality conditions for parabolic control problems have been studied, for instance, by Ahmed and Teo \cite{8.}, Cesari \cite{9.} and others for evolution and differential equations.

In the beginning of 1980s, the first papers on optimal control for variational inequalities appeared. In particular, such problems were intensively studied by Barbu \cite{10.} and Tiba \cite{11.}. Also, optimal control problems for some variational and hemivariational inequalities were considered in \cite{12.} and \cite{13.}, respectively. The optimal control problems for the subdifferential evolution inclusions have been examined in many works, see, e.g., \cite{14.}, \cite{15.}, \cite{16.} and references therein. More precisely, in \cite{15.} it was studied the Volterra integrodifferential evolution inclusions of nonconvolution type with time dependent unbounded operators and with both convex and nonconvex multivalued perturbations. In \cite{17.} the optimal control of parabolic variational inequalities is studied in the case where spatial domain is not necessarily bounded. In papers \cite{18.}, \cite{19.} the existence, uniqueness, and convergence of optimal control problems associated with parabolic variational inequalities of the second kind is studied. In \cite{20.}, the variational stability of optimal control problems involving subdifferential operators is studied. In \cite{21.} the optimal control problems for systems goverened by parabolic equations without initial conditions with controls in the coefficients were considered. The authors prove the existence of solutions to an optimal control Fourier problem for parabolic equations (without initial conditions), where the control functions occur in the coefficients of the state equation. They discuss the well-posedness of the problem and give a necessary optimality condition in the form of a generalized principle of Lagrange multipliers. Also, in \cite{22.}, the authors study optimal control problems governed by a class of semiliniear evolution equations including the so-called equations with memory. The evolutionary variational-hemivariational inequalities with applications to dynamic viscoelastic contact mechanics is considered in \cite{23.}. It is a complicated variational-hemivariational inequality of parabolic type with history-dependent operators which is studied here.

The most recent paper, known to us, in which the optimal control of history-dependent evolution inclusions is investigated is \cite{3.}. In this paper, a class of subdifferential evolution inclusions involving history-dependent operators was studied. The main idea of this paper is that the
existence and uniqueness results for such problems have been proved with removing the smallness condition. Also, the continuous dependence of the solution to these inclusions on the second member and initial condition were examined and the Bolza-type optimal control problem was studied. More precisely, the author at first consider the following problem
\begin{equation*}
\text{find}\quad w \in W \quad\text{such~that}
\end{equation*}
\begin{equation*}
w'(t)+A(t,w(t))+\partial\psi(t,w(t))\ni f(t)\quad\text{for~a.e}\,\,t\in(0,T),\quad w(0)=w_0.
\end{equation*}
Here, $\partial\psi$ denotes the Clarke generalized gradient of locally Lipschitz function $\psi(t,\cdot)$. It was shown that under some hypothesis the stated problem has a unique solution.

In the present paper, we turn our attention to the optimal control problem for a class of evolution inclusions with Volterra type operators. The mentioned above operator is supposed to be history-dependent. Furthermore, the stated problem is considered in the framework of an evolution triple of spaces. The existence and uniqueness results for the initial value problem of such inclusions have been proved in \cite{24.} As our main result, we establish the existence of a solution to the stated optimal control problem under some hypothesis on data-in. Based on the preceding papers, mentioned above, we can conclude that optimal control problems for subdifferential inclusions with Volterra type operators, considered here, have not been investigated yet. It serves us as one of the main motivations for the study of such kind of problems.

The outline of the paper is as follows. In Section 2, we formulated main notations and auxiliary facts. Statement of initial value problem for evolutionary subdifferential inclusions and the formulation of the results regarding the existence and uniqueness of its solutions is given in Section 3. The main results of this paper are stated in Section 4. In Section 5, we prove the main results.

\section{Preliminaries}
Let $T>0$ be an arbitrary fixed number. Let $V$ be a  separable reflexive Banach space  with the norm $\|\cdot\|$, $H$ be the Hilbert space  with the scalar product $(\cdot,\cdot)$ and the norm $|\cdot|$.
Assume that embedding $V\subset H$ is dense, continuous and compact.

By $V'$ and $H'$ we denote the  dual spaces to $V$ and $H$, respectively.
We assume, after an appropriate identification of functionals, that the space $H'$
is a subspace of $V'$.
Identifying $H'$ with $H$ by the Riesz-Fr\'{e}chet representation theorem, one usually writes
\begin{equation}\label{8I_vhv}
V\subset H\subset V'\,,
\end{equation}
where all embeddings are continuous and dense, and the first embedding is compact.
Note that in this case $\langle g,v\rangle_{V}=(g,v)$ for every
$v\in V,$ $g\in H,$ where $\langle\cdot,\cdot\rangle_{V}$  is the scalar product on the dual pair $\bigl[V',V\bigr]$.
Therefore, we will use notation $(\cdot,\cdot)$ instead of
$\langle\cdot,\cdot\rangle_{V}.$

Now, we introduce some functional spaces and spaces of distributions.

Let $X$ be an arbitrary Banach space with the norm $\|\cdot\|_X$.
By $C([0,T];X)$ we mean the Banach space of continuous functions $w:[0,T]\rightarrow X$ with the norm $\|w\|_{C(0,T;X)}:=\max\limits_{t\in[0,T]}\|w(t)\|_X$. By $L^q(0,T;X)$, where $q\geq1$, we denote the Banach space of measurable functions $w:(0,T)\rightarrow X$ such that $\|w(\cdot)\|_X\in L^q(0,T)$ with the norm $\|w\|_{L^q(0,T;X)}:=\Big(\int_0^T \|w(t)\|^q\,dt\Big)^{1/q}$. By $L^{\infty}(0,T;X)$ we denote the Banach space of measurable functions $w:(0,T)\rightarrow X$ such that $\|w(\cdot)\|_X\in L^{\infty}(0,T)$ with the norm $\|w\|_{L^{\infty}(0,T;X)}:= \operatorname*{ess~sup}\limits_{t\in(0, T)}\|w(t)\|_{X}$. By $D'(0,T;X)$ we mean the space of distributions on $D(0,T)$ with values on $X$, i.e., the space of  linear continuous (in weak topology on $X$) functionals on $D(0,T)$ with values on $X$ (hereafter $D(0,T)$ is space of test functions, that is,  the space of infinitely differentiable on $(0,T)$ functions with compact supports, equipped with the corresponding weak topology).

It is easy to see (by (\ref{8I_vhv})), that spaces $L^q(0,T;V)$, $L^{2}(0,T;H)$,  $L^{q'}(0,T;V')$, where $q>1$, $1/q+1/q'=1$, can be identified with the corresponding subspaces of $D'(0,T;V')$. In particular, this allows us to talk about derivatives $w'$ of functions $w$ from $L^q(0,T;V)$ and $L^{2}(0,T;H)$ in the sense
of distributions $D'(0,T;V')$  and
belonging of such derivatives to $ L^{q'}(0,T;V')$ or $L^{2}(0,T;H)$.

Let us define spaces
\begin{equation*}
H^1(0,T;H):=\{w \in L^{2}(0,T;H)\,\,\big|\,\, w'\in L^{2}(0,T;H)\},
\end{equation*}
\begin{equation*}
W_q^1(0,T;V):=\{w \in L^{q}(0,T;V)\,\,\big|\,\, w'\in L^{q'}(0,T;V')\},\quad q>1,\,\, q':=q/(q-1).
\end{equation*}

From the known results (see.,\cite[p.177-179]{25.})
it follows that $H^{1}(0,T;H)\subset C([0,T];H)$ and $W_q^1(0,T;V)\subset C([0,T];H)$.
Moreover, for every $w$ in $H^{1}(0,T;H)$ or $W_q^1(0,T;V)$
the function $t\mapsto |w(t)|^2$ is absolutely continuous
on $[0,T]$ and the following equality holds
\begin{equation}\label{8I_dy}
\frac{d}{dt}|w(t)|^2=2(w'(t),w(t)) \quad\text{for a.e.}\quad t\in [0,T].
\end{equation}

In this paper we use the following well-known facts.

\begin{proposition}[\cite{26.}, p.173]\label{proposition_1} Let $Y$ be a Banach space with the norm $\|\cdot\|_Y$, and $\{v_k\}_{k=1}^{\infty}$ be the sequence of elements of $Y$ which is weakly or $\ast$-weakly convergent to $v$ in $Y$. Then $\varliminf\limits_{k\to \infty}\|v_k\|_Y\geq \|v\|_Y$.
\end{proposition}

\begin{proposition}[\cite{27.}, Aubin Theorem]\label{proposition_2} Suppose that $q>1$, $r>1$ are some numbers, and $B_0, B_1,B_2$ are Banach spaces such that
$B_0\overset{c}{\subset}B_1\circlearrowleft B_2$
(here $\overset{c}{\subset}$ means compact embedding and $\circlearrowleft$ means continuous embedding). Then
\begin{equation*}
\{w\in L^q(0,T;B_0)\,|\,w'\in L^r(0,T;B_2)\}\overset{c}{\subset}\bigl(L^q(0,T;B_1)\cap
C([0,T];B_2)\bigr).
\end{equation*}
\end{proposition}
 We understand this embedding as follows: if the sequence $\{w_m\}_{m\in\mathbb{N}}$ is bounded in the space $L^q(0,T;B_0)$ and the sequence $\{w'_m\}_{m\in\mathbb{N}}$ is bounded in the space $L^r(0,T;B_2)$, then there exist a function $w\in C([0,T];B_2)\cap L^q(0,T;B_1)$ and a subsequence $w_{m_j}$ of the sequence $\{w_m\}$ such that $w_{m_j}\mathop{\longrightarrow}\limits_{j\to\infty}w$ in $C([0,T];B_2)$ and strongly in $L^q(0,T;B_1)$.

\section{Initial value problem for evolutionary subdifferential inclusions}

Let \label{8I_Functionallabel}$\Phi:V\rightarrow (-\infty, +\infty]$  be a proper functional, i.e., $\mathrm{dom}(\Phi):=\{v\in V\,|\,\ \Phi(v)<+\infty\}\neq\varnothing,$
which satisfies the following conditions
\bigskip
 \begin{description}
     \item[$(\mathcal{A}_1)$] \quad
  $ \Phi\bigr(\alpha v+(1-\alpha)w\bigl)\leq \alpha\Phi(v)+(1-\alpha)\Phi(w)$
   $ \forall\,  v, w\in V,\ \forall\, \alpha\in[0,1],$
   \bigskip
 i.e., the functional $\Phi$ is \emph{convex},
 \bigskip
 \item[$(\mathcal{A}_2)$]\quad
   $ v_{k}\mathop{\longrightarrow}\limits_{k\to\infty}v\,\,\,
   \textrm{ in }\,\,\, V\quad\Longrightarrow\quad
   \varliminf\limits_{k\to \infty}\Phi(v_{k})\geq \Phi(v),$

   \medskip
 i.e., the functional $\Phi$ is \emph{lower semi-continuous}.
\end{description}

\medskip
Let us recall (see, for example, \cite{1.}) that the \emph{subdifferential} of functional $\Phi$ is a map
$\partial\Phi: V\rightarrow 2^{{V}'}$, defined as follows
 \begin{equation*}
   \partial\Phi(v):=\{v^{\ast}\in{V}'\ | \
   \Phi(w)\geq\Phi(v)+(v^{\ast}, w-v)\quad\forall\ w\in {V}\},\quad v\in {V},
 \end{equation*}
 and the  \emph{domain} of the subdifferential $\partial\Phi$ is the set
 $D(\partial\Phi):=\{v\in{V}\,|\, \partial\Phi(v)\neq\varnothing\}$.
 We identify the subdifferential $\partial\Phi$ with its graph, assuming that
 $[v,v^{\ast}]\in\partial\Phi$ if and only if $v^{\ast}\in\partial\Phi(v)$,
 i.e., $\partial\Phi=\{[v,v^*] \ | \ v\in D(\partial\Phi),\ v^*\in \partial\Phi(v))\}$.
 R. Rockafellar (see \cite[theorem~A]{28.}) proves that
 the subdifferential $\partial\Phi$ is a \emph{maximal monotone operator}, that is,
 \begin{equation}\label{monoton}
    (v^{\ast}_{1}-v^{\ast}_{2},v_{1}-v_{2})\geq0\qquad \forall\
    [v_{1}, v^{\ast}_{1}],\ [v_{2}, v^{\ast}_{2}]\in\partial\Phi,
 \end{equation}
  and for every element $[v_{1}, v_{1}^{\ast}]\in{V}\times{V}'$ we have the following implication
 \begin{equation*}
    (v^{\ast}_{1}-v^{\ast}_{2}, v_{1}-v_{2})\geq0\quad
    \forall\ [v_{2}, v^{\ast}_{2}]\in\partial\Phi\ \ \Rightarrow\ \ [v_{1}, v^{\ast}_{1}]\in\partial\Phi.
 \end{equation*}

\medskip
Additionally, assume that the following condition holds

\medskip\noindent
$(\mathcal{A}_3)$~     there exist  constants $p\geq2,\ K>0$ such that
     \begin{equation*}
       \Phi(v)\geq K\|v\|^{p}\quad\forall\ v\in\mathrm {dom}(\Phi);
     \end{equation*}
     moreover, $\Phi(0)=0.$

\medskip
\begin{remark}\label{remark_2}
Condition $(\mathcal{A}_3)$ implies that $\Phi(v)\geq\Phi(0)+(0,v-0)\,\,\,\,\forall\,v\in V$, hence $[0,0]\in\partial\Phi$.
\end{remark}

\medskip

Suppose that

\medskip
$(\mathcal{A}_4)$~ $\mathcal{B}: C([0,T];H)\rightarrow L^\infty(0,T;H)$ is a bounded linear operator such that,
 for a.e. $t\in(0,T)$ and for any $w_1,w_2\in C([0,T];H)$, the following inequality holds
\begin{equation*}
|\mathcal{B}(w_1)(t)-\mathcal{B}(w_2)(t)|\leq L\int\limits_0^t |w_1(s)-w_2(s)|\,ds,
\end{equation*}
 where $L\geq0$ is some a constant;  moreover, $\mathcal{B}(0)=0$.

It is easy to varify that condition $(\mathcal{A}_4)$  implies
\begin{equation}\label{operator_B}
|\mathcal{B}(w)(t)|\leq L\int\limits_0^t|w(s)|\,ds
\end{equation}
for a.e. $t\in(0,T)$ and for any $w\in C([0,T];H)$.

Operator $\mathcal{B}$ is called a {\it Volterra} type operator.

\medskip
\begin{remark}\label{remark_3}
An example of operator $\mathcal{B}$ in $(\mathcal{A}_4)$ is the operator $\widehat{\mathcal{B}}$, introduced in Section 1. Another important example of operator $\mathcal{B}$ is the following operator
\begin{equation}\label{example}
\widetilde{\mathcal{B}}(w)(t):=\widetilde{B}\big(t,\int_0^t \widetilde b(t,s,w(s))\,ds\big),\quad t\in(0,T),
\end{equation}
where $\widetilde{B}:(0,T)\times H\rightarrow H$, $\widetilde b:\Pi\times H\rightarrow H$ are maps which satisfy the following conditions
\begin{itemize}
\item[$1)$]
for any $v\in H$ the map $\widetilde{B}(\cdot,v):(0,T)\rightarrow H$ is measurable, and
there exists a constant $L_1\geq0$ such that the following inequality holds
\begin{equation}\label{condition_2}
|\widetilde{B}(t,v_1)-\widetilde{B}(t,v_2)|\leq L_1|v_1-v_2|
\end{equation}
\medskip
for a.e. $t\in(0,T)$ and for all $v_1,v_2\in H$; in addition, $\widetilde{B}(t,0)=0$ for a.e. $t\in(0,T)$.
\item[$2)$]
for any $v\in H$ the map $\widetilde b(\cdot,\cdot,v):\Pi\rightarrow H$ is measurable, and there exists a constant $L_2\geq0$ such that the following inequality holds
\begin{equation}\label{condition_1}
|\,\widetilde b(t,s,v_1)-\widetilde b(t,s,v_2)|\leq L_2|v_1-v_2|
\end{equation}
for a.e. $(t,s)\in\Pi$ and for all $v_1,v_2\in H$; in addition, $\widetilde b(t,s,0)=0$ for a.e. $(t,s)\in\Pi$.
\end{itemize}

Let us show that operator $\widetilde{\mathcal{B}}$, defined in (\ref{example}), satisfies condition $(\mathcal{A}_4)$ with $L=L_1L_2$. Indeed, based on  $(\ref{condition_2})$  and $(\ref{condition_1})$ we have
\begin{equation*}
|\widetilde{\mathcal{B}}(w_1)(t)-\widetilde{\mathcal{B}}(w_2)(t)|\leq\big|\widetilde{B}\big(t,\int_0^t \widetilde b(t,s,w_1(s))\,ds\big)-
\widetilde{B}\big(t,\int_0^t \widetilde b(t,s,w_2(s))\,ds\big)\big|\leq$$ $$\leq
 L_1\Big|\int_0^t \widetilde b(t,s,w_1(s))\,ds-\int_0^t \widetilde b(t,s,w_2(s))\,ds\Big|\leq L_1L_2\int_0^t|w_1(s)-w_2(s)|\,ds.
\end{equation*}
\end{remark}

\medskip

Let us consider {\it\textbf{evolutionary subdifferential inclusion}}
\begin{equation}\label{Equation14}
  y'(t)+\partial\Phi\bigl(y(t)\bigr)+\mathcal{B}(y)(t)\ni f(t),\quad t\in (0,T),
\end{equation}
where $f: (0,T)\to V'$ is a given measurable function, $y:[0,T]\rightarrow V$ is an unknown function.

\medskip
Let conditions $(\mathcal{A}_1)$ -- $(\mathcal{A}_4)$ hold, and $f\in L^{p'}(0,T;V')$, where $p'=p/(p-1)$.

\medskip
\begin{definition}\label{def1}
\emph{The function $y$ is a} solution \emph{of variational inequality \eqref{Equation14} if it satisfies the following conditions}
   \begin{itemize}
   \item[\textbf{1)}]\label{8I_1)}
     $y\in W_p^1(0,T;V)$ \ (\emph{then} $y\in C([0,T];H)$);
   \item[\textbf{2)}]\label{8I_2)}
     $y(t)\in D(\partial\Phi)$ \emph{for a.e.} $t\in (0,T)$;
   \item[\textbf{3)}]\label{8I_3)}
     \emph{there exists a function} $g \in L^{p'}(0,T;V')$ \emph{such that for a.e.} $t\in (0,T)$ \emph{we have}
     $g(t)\in \partial\Phi\bigl(y(t)\bigr)$ \emph{and}
     \begin{equation}\label{Equation15}
       y'(t)+g(t)+\mathcal{B}(y)(t)=f(t)\quad\textrm{\emph{in}}\quad V'.
     \end{equation}
   \end{itemize}
\end{definition}
By~$\textbf{P}(\Phi,\mathcal{B},f,y_0)$ we denote the {\it\textbf{problem}} of finding a solution $y$ of the variational inequality
  \eqref{Equation14} that satisfies the following condition
\begin{equation}\label{Equation16}
y(0)=y_0,
\end{equation}
where $y_0\in H $ is given.

\medskip

\begin{remark}\label{remark_4} The problem~$\textbf{P}(\Phi,\mathcal{B},f,y_0)$ can be replaced by the following one (see Introduction). Let $K$ be a convex, closed set in $V$, $A:V\to V'$ be a
monotone, bounded and semi-continuous operator such that
$ (A(v),v)\geq \widetilde{K}_1\|v\|^p\quad  \forall v \in V$, where $\widetilde{K}_1=$ const $>0$, and $f\in L^{p'}(0,T;V'),\ y_0\in H$.
The problem is to find a function $y \in W_p^1(0,T;V)$ such that $y(0)=y_0$ and, for a.e. $t \in (0,T)$, we have $y(t)\in K$ and
\begin{equation*}
(y'(t)+A(y(t))+\mathcal{B}(y)(t), v-y(t))\geq(f(t),v-y(t))\quad \forall\,\, v\in K.
\end{equation*}
\end{remark}

In \cite{24.} we obtain the following results.

\begin{theorem}[\cite{24.}, Theorem 1]\label{th1} {\sl Let conditions $(\mathcal{A}_1)$ -- $(\mathcal{A}_4)$ hold, and $f\in L^{p'}(0,T;V')$, $y_0\in H$. Then the problem~$\textbf{P}(\Phi,\mathcal{B},f,y_0)$ has no more than one solution}.
\end{theorem}

\begin{theorem}[\cite{24.}, Theorem 2]\label{th2}  Let conditions
 $(\mathcal{A}_1)$--$(\mathcal{A}_4)$ hold, and $f\in L^{2}(0,T;H)$, $y_0\in \mathrm{dom}(\Phi)$.
  Then the problem~$\textbf{P}(\Phi,\mathcal{B},f,y_0)$ has a unique solution which belongs to $L^\infty(0,T;V)\cap H^{1}(0,T;H)$,
  and satisfies the following estimate
  \begin{equation}\label{Equation17}
     \max_{t\in [0,T]}|y(t)|^{2}+  \mathop{{\mathrm {ess}}\sup}\limits_{t\in [0,T]}\|y(t)\|^p+\int_{0}^{T}|y'(t)|^{2}\,dt\leq$$ $$\leq C_1\Bigl(|y_0|^2+\Phi(y_0)+\int_{0}^{T}|f(t)|^{2}\,dt\Bigr),
   \end{equation}
  where $C_1$ is a positive constant, which depends only on $K,L$, and $T$.
\end{theorem}

\section{Statement of the optimal control problem and the main result}

Let $H^{\ast}$ be a Hilbert space with the scalar product $(\cdot,\cdot)_{H^{\ast}}$ and the corresponding norm $\|\cdot\|_{H^{\ast}}:=\sqrt{(\cdot,\cdot)_{H^{\ast}}}$. Let us consider the space of controls
\begin{equation*}
U:=\big\{u\in L^2_{\textrm{loc}}(0,T;H^{\ast})\ \big| \ \int_0^T\omega(t)\|u(t)\|^2_{H^{\ast}}dt<\infty\big\},
\end{equation*}
 where $\omega\in C(0,T),$ and $\omega(t)>0$ for all $t\in (0,T)$.  It is a Hilbert space with the scalar product and norm, respectively,
\begin{equation*}
(u_1,u_2)_{U}:=\int\limits_0^T\omega(t)(u_1(t),u_2(t))_{H^{\ast}}\,dt,\quad
\|u\|_U:=\Bigl(\int\limits_0^T\omega(t)\|u(t)\|^2_{H^{\ast}}\,dt\Bigr)^{1/2},
\end{equation*}
for all $u_1,\,u_2,\,u\in U$.

\medskip
Let $U_{\partial}\subset U$ be a convex closed set, which is called the set of admissible controls.

Further, we assume that the $\Phi$, $\mathcal{B}$ satisfy conditions $(\mathcal{A}_1)$ -- $(\mathcal{A}_4)$, and

\medskip\noindent
 $(\mathcal{C})$~   $C:U\to L^2(0,T;H)$ is a linear continuous operator

\medskip\noindent
     $(\mathcal{F})$~     $f\in L^2(0,T;H),\quad y_0\in \mathrm{dom}(\Phi)$.

\medskip
For a given control $u\in U_{\partial}$ the state $y(t),\, t\in [0,T],$ (which can be also denoted by $y(u)$ or $y(\cdot;u)$) of the control evolutionary system  is described by the solution of the evolutionary variational inequality
\begin{equation}\label{Equation21}
y'(t)+\partial\Phi\bigl(y(t)\bigr)+\mathcal{B}(y)(t)\ni f(t)+Cu(t),\quad t\in (0,T),
\end{equation}
with initial condition
\begin{equation}\label{pochumovazad}
y(0)=y_0,
\end{equation}
that is, $y$ is a solution to the problem~$\textbf{P}(\Phi,\mathcal{B},f+Cu,y_0)$.

From  Theorem \ref{th1} and \ref{th2} it follows that there exists a unique solution $y\in L^\infty(0,T;V)\cap H^{1}(0,T;H)$ to the problem~$\textbf{P}(\Phi,\mathcal{B},f+Cu,y_0)$ which satisfies the following estimate
  \begin{equation}\label{Equation17}
      \max_{t\in [0,T]}|y(t)|^{2}+ \mathop{{\mathrm {ess}}\sup}\limits_{t\in [0,T]}\|y(t)\|^p+\int_{0}^{T}|y'(t)|^{2}\,dt\leq$$ $$\leq
           C_1\Bigl(|y_0|^2+\Phi(y_0)+2\int_{0}^{T}\big(|f(t)|^{2}+|Cu|^{2}\big)\,dt\Bigr),
   \end{equation}
  where $C_1$ is a positive constant, which depends only on $K,L$, and $T$.

\medskip
Let the following condition holds
\begin{description}
     \item[$(\mathcal{G})$] \label{8I_(Chpt05:0,Tect01:I)}~
     $G: C([0,T];H)\rightarrow\mathbb{R}$ be a lower semi-continuous, bounded bellow functional, i.e.,
     $$\inf_{z\in C([0,T];H)} G(z)>-\infty.$$
\end{description}

Let us define the cost functional $J:U\rightarrow\mathbb{R}$ by the following rule
\begin{equation}\label{Equation22}
J(u):=G(y(u))+\mu\|u\|^2_U,\quad u\in U,
\end{equation}
where $\mu>0$ is a constant, $y(u)$ is the solution to the problem~$\textbf{P}(\Phi,\mathcal{B},f+Cu,y_0)$.

The \textbf{optimal control problem} is to find $u^{\ast}\in U_{\partial}$ such that
\begin{equation}\label{Equation23}
J(u^{\ast})=\inf_{u\in U_{\partial}}J(u).
\end{equation}
Later, this problem will be called as a problem \eqref{Equation23}.

\medskip
The main result of this paper is stated in the following theorem.

\begin{theorem}\label{th3} Let conditions $(\mathcal{A}_1)$--$(\mathcal{A}_4)$, $(\mathcal{C})$, $(\mathcal{F})$ and $(\mathcal{G})$ hold. Then the problem $(\ref{Equation23})$ has a solution.
\end{theorem}

The proof of this theorem will be presented in the next section.

\section{Proof of the main result}

Let us prove Theorem \ref{th3}. Since the cost function $J$ is bounded below, it implies that there exists the sequence $\{u_k\}_{k=1}^{\infty}\in U_{\partial}$ such that
\begin{equation}\label{Equation635}
J(u_k)\mathop{\longrightarrow}\limits_{k\to\infty}{}J_{\ast}:=\inf_{u\in U_{\partial}}J(u)>-\infty
\end{equation}
\quad
Thus it means that the sequence $\bigl\{J(u_k)\bigr\}_{k=1}^{\infty}$ is bounded. Taking into account (\ref{Equation22}), one can obtain that the sequence $\{u_k\}_{k=1}^{\infty}$ is bounded in $U$, i.e.,
\begin{equation}\label{Equation63}
\|u_k\|_U\leq C_2,\quad k\in\mathbb{N},
\end{equation}
where $C_2$ is a constant, which does not depend on $k$.

Since $C:U\to L^2(0,T;H))$ is linear continuous operator, then $\bigl\{Cu_k\bigr\}_{k=1}^{\infty}$ is bounded in $L^2(0,T;H)$, that is
\begin{equation}\label{Equation603}
\|Cu_k\|_{L^2(0,T;H)}\leq C_3,\quad k\in\mathbb{N}.
\end{equation}
where $C_3$ is a constant, which does not depend on $k$.

For each $k\in \mathbb{N}$, denote $y_k:=y(u_k)$, i.e., $y_k$ is a solution to the problem $\textbf{P}(\Phi,\mathcal{B},f+Cu_k,y_0)$. Taking into account condition $(\mathcal{C})$ and $(\mathcal{F})$, from Theorem \ref{th2} it follows that for each $k\in \mathbb{N}$, $y_k\in L^{\infty}(0,T;V)\cap H^1(0,T;H)$,  for a.e. $t\in (0,T)$,  $y_k(t)\in D(\partial\Phi)$ and
\begin{equation}\label{Equation65}
y_k'(t)+\partial\Phi(y_k(t))+\mathcal{B}(y_k)(t)\ni f(t)+Cu_k(t)\quad \textrm{in}\ H,
\end{equation}
\begin{equation}\label{Equation605}
y_k(0)=y_0.
\end{equation}
Also, the following estimate holds
\begin{equation}\label{Equation67}
\max_{t\in [0,T]}|y_k(t)|^{2}+\mathop{{\mathrm {ess}}\sup}\limits_{t\in [0,T]}\|y_k(t)\|^p+\int_{0}^{T}|y'_k(t)|^{2}\,dt\leq
$$ $$\leq C_1\Bigl(|y_0|^2+\Phi(y_0)+2\int_{0}^{T}\bigl(|f(t)|^{2}+|Cu_k(t)|^{2}\bigr)\,dt\Bigr),
\end{equation}
where $C_{1}$ is a positive constant which  depends only on $K,L$, and $T$.

Moreover, from the proof of Theorem \ref{th2} it follows that there exists the sequence $g_k\in L^{2}(0,T;H)$ such that, for a.e. $t\in(0,T)$, $g_k(t)\in\partial\Phi\bigl(y_k(t)\bigr)$  and
\begin{equation}\label{Equation66}
y'_k(t)+g_k(t)+\mathcal{B}(y_k)(t)=f(t)+Cu_k(t)\,\,\,\text{in}\,\,H.
\end{equation}

From (\ref{Equation603}), (\ref{Equation67}),  and $(\mathcal{F})$ it follows that
\begin{equation}\label{Equation69}
\{y_k\}_{k=1}^{\infty}\,\,\,\text{is bounded in}\,\,L^{\infty}(0,T;V)\cap C([0,T];H),
\end{equation}
\begin{equation}\label{Equation70}
\{y'_k\}_{k=1}^{\infty}\,\,\,\text{is bounded in}\,\,L^{2}(0,T;H).
\end{equation}

Taking into account (\ref{operator_B}) and (\ref{Equation69}), we obtain
\begin{equation}\label{Equation71}
\{\mathcal{B}(y_k)\}_{k=1}^{\infty}\,\,\,\text{is bounded in}\,\,L^{2}(0,T;H).
\end{equation}
\quad
From (\ref{Equation66}), taking into account (\ref{Equation603}), (\ref{Equation70}), (\ref{Equation71}), and $(\mathcal{F})$, we obtain
\begin{equation}\label{Equation72}
\{g_k\}_{k=1}^{\infty}\,\,\,\text{is bounded in}\,\,L^{2}(0,T;H).
\end{equation}
\quad
Let us recall that spaces $V$ and $H$ are reflexive. So, from (\ref{Equation63}), (\ref{Equation69}), (\ref{Equation70}), (\ref{Equation72}), taking into account
Proposition \ref{proposition_2}, it follows that there exists a subsequence of a sequence $\{(u_k, y_k, g_k)\}_{k=1}^{\infty}$ (which we denote by $\{(u_k, y_k, g_k)\}_{k=1}^{\infty}$) and functions $u_{\ast}\in U_{\partial}$, $y\in L^{\infty}(0,T;V)\cap H^1(0,T;H)$ (then $y\in C([0,T];H)$), and $g\in L^2(0,T;H)$ such that
\begin{equation}\label{Equation73}
u_k\mathop{\longrightarrow}\limits_{k\to\infty}u_{\ast}\,\,\text{weakly in}\ \ U,
\end{equation}
\begin{equation}\label{Equation74}
y_k\mathop{\longrightarrow}\limits_{k\to\infty}y\,\,\ast\text{-weakly in}\,\,L^{\infty}(0,T;V),\,\,\,\text{weakly in}\,\,H^1(0,T;H),
\end{equation}
\begin{equation}\label{Equation75}
y_k\mathop{\longrightarrow}\limits_{k\to\infty}y\quad\text{in}\ \ C([0,T];H),
\end{equation}
\begin{equation}\label{Equation77}
g_k\mathop{\longrightarrow}\limits_{k\to\infty}g\quad\text{weakly in}\ \ L^2(0,T;H).
\end{equation}

Based on $(\mathcal{A}_4)$ and \eqref{Equation75}, we obtain
\begin{equation*}
\mathop{\textrm{ess sup}}\limits_{t\in[0,T]}|\mathcal{B}(y_k)(t)-\mathcal{B}(y)(t)|\leq L\int_0^T |y_k(s)-y(s)|\,ds\mathop{\longrightarrow}\limits_{k\to\infty} 0,
\end{equation*}
that is
\begin{equation}\label{Equation76}
\mathcal{B}(y_k)\mathop{\longrightarrow}\limits_{k\to\infty} \mathcal{B}(y)\quad
\text{strongly in}\ \ L^\infty(0,T;H).
\end{equation}

From $(\mathcal{C})$ and (\ref{Equation73}) it follows
\begin{equation}\label{Equation765}
Cu_k\mathop{\longrightarrow}\limits_{k\to\infty} Cu_{\ast}\quad
\text{weakly in}\ \ L^2(0,T;H).
\end{equation}

Let $v\in H, \varphi \in C([0,T])$ be arbitrary functions. For a.e. $t\in (0,T)$ we multiply equality \eqref{Equation66} by $v$,  and
then we multiply the obtained equality by
$\varphi$ and integrate it with respect to $t\in[0,T]$. As a result,
 we obtain
\begin{align}\label{step3_1}
\int_0^T(y'_k(t),v\varphi(t))\,dt+\int_0^T(g_k(t),v\varphi(t))\,dt
+\int_0^T\!\big(\mathcal{B}(y_k)(t),v\varphi(t)\big)\,dt
\notag \\
=\int_0^T(f(t)+Cu_k(t),v\varphi(t))\,dt,\quad k\in \mathbb{N}.
\end{align}
\quad
Taking into account (\ref{Equation73})  -- (\ref{Equation765}), we pass to the limit in (\ref{Equation66}) as $k\rightarrow\infty$. As a result, since $v\in H, \varphi \in C([0,T])$ are arbitrary, for a.e. $t\in (0,T)$ we obtain
\begin{equation*}
y'(t)+g(t)+\mathcal{B}(y)(t)=f(t)+(Cu_{\ast})(t)\,\,\,\text{in}\ \ H.
\end{equation*}

In order to complete the proof of the theorem  it remains to show that
$y(t)\in D(\partial\Phi)$ and  $g(t)\in\partial\Phi\bigl(y(t)\bigr)$ for a.e. $t\in (0,T)$.

Let $k\in \mathbb{N}$ be an arbitrary number.
Since $y_k(t)\in D(\partial \Phi)$ and $g_{k}(t)\in\partial\Phi\bigl(y_{k}(t)\bigr)$ for a.e. $t\in (0,T)$,
we obtain, using the monotonicity of the subdifferential  $\partial\Phi$, that for a.e. $t\in (0,T)$ the following equality holds
  \begin{equation}\label{inequality139}
    (g_{k}(t)-v^{\ast},y_{k}(t)-v)\geq0\quad\forall\, [v, v^{\ast}]\in\partial\Phi.
 \end{equation}
 \quad
 Let $\sigma \in (0,T)$, $h>0$ be arbitrary numbers such that $\sigma-h\in (0,T)$.
 We integrate  \eqref{inequality139} on $(\sigma-h,\sigma)$ and obtain
 \begin{equation}\label{Chpt05:Sect01:Equ41}
    \int_{\sigma-h}^{\sigma}(g_{k}(t)-v^{\ast},y_{k}(t)-v)\,dt \geq0
    \quad\forall\, [v, v^{\ast}]\in\partial\Phi.
 \end{equation}
 \quad
 Now, according to \eqref{Equation74} -- \eqref{Equation77}, we pass to the limit in \eqref{Chpt05:Sect01:Equ41} as $k\to\infty$.
 As a result we obtain
  \begin{align}\label{Chpt05:Sect01:Equ42}
       \int_{\sigma-h}^{\sigma}(g(t)-v^{\ast},y(t)-v)\,dt\geq0
   \quad\forall\, [v, v^{\ast}]\in\partial\Phi.
 \end{align}

The monograph \cite[Theorem~2, p.~192]{26.} and \eqref{Chpt05:Sect01:Equ42}
 imply that for every $[v, v^{\ast}]\in\partial\Phi$ there exists a set
 $R_{[v,v_*]}\subset (0,T)$  of measure zero such that for all $\sigma\in (0,T)\setminus R_{[v,v_*]}$ we have
$y(\sigma)\in V,\ g(\sigma)\in H$ and
 \begin{equation}\label{1468'}
 0\leq  \lim\limits_{h\to+0}\frac{1}{h}\int_{\sigma-h}^{\sigma}\big(g(t)-v^*,y(t)-v\big)\,dt=
  \big(g(\sigma)-v^*,y(\sigma)-v\big).
   \end{equation}
   \quad
Let us show that there exists a set of measure zero  $R\subset (0,T)$ such that
 \begin{equation}\label{1468}
 \forall \sigma\in (0,T)\setminus R: \qquad
    \big(g(\sigma)-v^*,y(\sigma)-v\big) \geq 0 \quad \forall\, [v,v^{\ast}]\in\partial\Phi.
   \end{equation}
   \quad
Since  $V$ and $V'$ are separable spaces, there exists a countable set $F\subset \partial\Phi$
which is dense in $\partial\Phi$.
Let us denote by $R:=\mathop{\cup}\limits_{[v,v^*]\in F}R_{[v,v_*]}$.
Since the set $F$ is countable, and any countable union of sets of measure zero is a set of measure zero,
then $R$ is a set  of measure zero.
Therefore, for any $\sigma\in S\setminus R$ inequality
\begin{equation*}
\big(g(\sigma)-v^*,y(\sigma)-v\big) \geq 0\quad \forall [v, v^{\ast}]\in F
\end{equation*}
 holds.
Let $[\widehat{v},\widehat{v}^*]$ be an arbitrary element from $\partial\Phi$.
Then, since $F$ is dense in $\partial\Phi$, we have
the existence of a sequence $\{[v_l,v_l^*]\}_{l=1}^\infty$ such that
$v_l\to \widehat{v}$ in $V$, $v_l^*\to \widehat{v}^*$ in $V'$, and
\begin{equation}\label{inequality142}
\forall \sigma\in (0,T)\setminus R: \qquad
  (g(\sigma)-v_l^{\ast},y(\sigma)-v_l)\geq0\quad\forall\, l\in \mathbb{N}.
 \end{equation}
 \quad
Thus, passing to the limit in inequality \eqref{inequality142} as $l\to \infty$, we get
$(g(\sigma)-\widehat{v}^*,y(\sigma)-\widehat{v})\geq0$ $\forall \sigma\in (0,T)\setminus R$.
Therefore,  inequality \eqref{1468} holds.
From this, according to the maximal monotonicity of $\partial\Phi$, we obtain
that $[y(t), g(t)]\in\partial\Phi$ for a.e. $t\in (0,T)$, that is
 $y(t)\in D(\partial\Phi)$ and $g(t)\in \partial\Phi(y(t))$ for a.e. $t\in(0,T)$. Thus, the function $y$ is a solution to the problem $\textbf{P}(\Phi,\mathcal{B},f+Cu_{\ast},y_0)$, i.e.,  $y=y(u_*)$ is the state
of the control evolutionary system for a given control $u_*$.

It remains to show that $u_{\ast}$ is a minimizing element of the functional $J$. Indeed, the functional $G$ is lower semi-continuous in $C([0,T];H)$, then (\ref{Equation75}) implies that
\begin{equation}\label{Equation78}
\varliminf\limits_{k\to \infty}G(y(u_k))\geq G(y(u_*))\,.
\end{equation}
\quad
Also, (\ref{Equation73}) and Proposition \ref{proposition_1}  yield
\begin{equation}\label{Equation79}
\varliminf\limits_{k\to \infty}\|u_k\|_U\geq\|u_{\ast}\|_U\,.
\end{equation}
\quad
From (\ref{Equation22}), (\ref{Equation635}), (\ref{Equation78}), and (\ref{Equation79}) we obtain that
\begin{align*}
\inf_{u\in U_{\partial}}J(u)=\lim_{k\to \infty}J(u_k)\geq\varliminf\limits_{k\to \infty}G(y(u_k))+\mu\varliminf\limits_{k\to \infty}\|u_k\|_U^2\geq\\
\geq G(y(u_*))+\mu\|u_*\|_U^2= J(u_{\ast}).
\end{align*}
\quad
Thus, we obtain that $u_{\ast}$ is a solution to the problem (\ref{Equation23}).
Hence, Theorem \ref{th3} is proved.

\section{Conclusions}

Our aim in the present paper was to study  the existence of a solution of the optimal control problem for a class of evolution inclusions with Volterra type operators. We recall that the Volterra type operator is assumed to be history-dependent. The main goal was achieved by satisfying necessary assumptions on input data. The result of our paper extends results studied in \cite{24.}, where the existence and uniqueness of the weak solution for the initial value problem of evolution inclusions with Volterra type operators was proved.

\begin{acknowledgements}
We wish to thank the referees and the associate editor for their helpful comments and suggestions.
\end{acknowledgements}

\end{document}